# CONDITIONAL DENSITY ESTIMATION IN A REGRESSION SETTING[1]

BY SAM EFROMOVICH

*The University of Texas at Dallas*

Regression problems are traditionally analyzed via univariate characteristics like the regression function, scale function and marginal density of regression errors. These characteristics are useful and informative whenever the association between the predictor and the response is relatively simple. More detailed information about the association can be provided by the conditional density of the response given the predictor. For the first time in the literature, this article develops the theory of minimax estimation of the conditional density for regression settings with fixed and random designs of predictors, bounded and unbounded responses and a vast set of anisotropic classes of conditional densities. The study of fixed design regression is of special interest and novelty because the known literature is devoted to the case of random predictors. For the aforementioned models, the paper suggests a universal adaptive estimator which (i) matches performance of an oracle that knows both an underlying model and an estimated conditional density; (ii) is sharp minimax over a vast class of anisotropic conditional densities; (iii) is at least rate minimax when the response is independent of the predictor and thus a bivariate conditional density becomes a univariate density; (iv) is adaptive to an underlying design (fixed or random) of predictors.

**1. Introduction.** Let $(Y_l, X_l)$, $l = 1, \ldots, n$, be independent pairs of observations (bivariate data). We would like to analyze a relationship (association) between variables $X_l$ (the predictor) and $Y_l$ (the response) that allows one to quantify the input of $X_l$ on $Y_l$. To simplify the problem, the nonparametric regression literature recommends analysis of the association via the conditional expectation of the response given the predictor because this

Received November 2005; revised January 2007.
[1]Supported in part by BIFAR and by NSF Grants DMS-02-43606 and DMS-06-04558.

implies estimation of a well-understood univariate function. In practical applications, this simplification may or may not fully describe the association; see discussion in [1, 5, 14, 16, 21, 31]. In general, the conditional density of the response given the predictor describes the ultimate association between the response and the predictor. However, this is a bivariate function and its estimation is complicated by the curse of dimensionality, the latter necessitating development of optimal estimators. The literature on such estimators is next to nothing, and the aim of this article is to develop minimax and oracle theory of estimation of conditional densities.

Let us formulate the problem of estimation of the conditional density (in what follows, the abbreviation c.d. will often be used) considered in this paper. We would like to estimate the c.d. of the response $Y$ given the predictor $X$ in the following regression settings. First, we need to take into account two possible models of design of predictors. The first model is where pairs of observations are independent samples from a pair of two random variables $Y$ and $X$. Then, if the joint density $f(y,x)$ exists and the marginal density $p(x) := \int_{-\infty}^{\infty} f(y,x)\,dy$ of the predictor is positive, we are estimating the conditional density

$$(1.1) \qquad f(y|x) := \frac{f(y,x)}{p(x)}.$$

It is traditional to refer to this design as random and to the marginal density $p$ as the design density, regardless of the fact that it may be known or unknown to the statistician. The second model is where predictors are created by a deterministic procedure and then responses are generated according to a conditional density $f(y|x)$. This is the case of a so-called *fixed design*. A discussion of these two designs can be found in [5, 31]; an interesting probabilistic point of view is presented in [1].

We also need to take into account that (i) the response can be either bounded or unbounded (the former case is typical in practical applications and the latter is of theoretical interest); (ii) the smoothness of the c.d. $f(y|x)$ may depend on the direction (it can be anisotropic), and moreover, if the response and predictor are independent, then $f(y|x) = f(y)$; (iii) different losses can be used to evaluate the quality of estimation of the c.d. All these issues will be explored in this paper.

The level of known results on c.d. estimation is not on a par with the theory of multivariate density estimation. The latter is the reason why using (1.1) has been the main approach to assess the optimality of a c.d. estimator. To give an example of how this formula is used in the literature, let us note that an isotropic bivariate density with two derivatives for each component can be estimated with Mean Integrated Squared Error (MISE) of order $n^{-2/3}$ and then if the design density is sufficiently smooth (say it is twice differentiable), this implies that the conditional density can also be evaluated with MISE of order $n^{-2/3}$. While such an approach is legitimate, it has obvious limitations. In particular, it cannot resolve many basic issues like how smoothness of the design density affects estimation of the conditional density or how to consider a classical fixed design regression.

Formula (1.1) has also been an inspiration for creating ad hoc estimators of the c.d. with the main theoretical emphasis on the bias-variance analysis. The interested reader can find a historical overview of this and related approaches in the books [5, 15, 34]; other relevant references are [2, 16, 17, 21, 22, 23, 26, 27, 28, 35].

The content of the article is as follows. Section 2 presents the setting. Sections 3 and 4 describe new sharp minimax lower bounds under $L_2([0,1]^2)$ and $L_2((-\infty,\infty) \times [0,1])$ losses, respectively. A c.d. estimator is defined in Section 5. An oracle inequality which shows how well the estimator matches an oracle that knows an estimated c.d. is presented in Section 6. Minimax properties of the estimator are established in Section 7. Optimal design of predictors for controlled experiments is explored in Section 8. Discussion of the results obtained, including analysis of real datasets, can be found in Section 9. Proofs are deferred to Section 10.

The following notation will be used throughout the article: $i$ always denotes the complex unit, that is, $i^2 = -1$; $\text{Re}\{\cdot\}$ is the real part; $o(1)$'s are generic sequences in $n$ such that $o(1) \to 0$ as $n \to \infty$; $Q$ is a positive constant; $C$'s are generic positive constants; $(x)_+ := \max(0, x)$; $\lfloor x \rfloor$ is the integer part of $x$; $I(\cdot)$ is the indicator; the cosine basis on $[0,1]$ is denoted by $\varphi_0(x) := 1, \varphi_j := 2^{1/2} \cos(\pi j x)$, $j = 1, 2 \ldots$. We shall use two different loss functions to study the performance of an estimator $\tilde{f}(y|x)$: $L_2([0,1]^2)$ loss, which is $\int_{[0,1]^2} (\tilde{f}(y|x) - f(y|x))^2 \, dy \, dx$, and $L_2((-\infty,\infty) \times [0,1])$ loss, which is $\int_0^1 [\int_{-\infty}^{\infty} (\tilde{f}(y|x) - f(y|x))^2 \, dy] \, dx$. If these two loss functions are considered simultaneously, then the area of integration is not written with the understanding that it corresponds to an underlying loss.

**2. Considered model.** Observations are $n$ pairs $\{(Y_l, X_l), l = 1, \ldots, n\}$ which are generated according to one of the following two designs. (i) *Random design*. The pairs are independent samples from a pair $(Y, X)$ of two random variables (the response and the predictor) with the joint density $f(y, x)$. Set $p(x) := \int_{-\infty}^{\infty} f(y, x) \, dy$ for the marginal (design) density of the predictor $X$. Assume that $p(x)$ is positive over its support. Then the problem is to find a corresponding conditional density (c.d.) $f(y|x) := f(y, x)/p(x)$. (ii) *Fixed design*. Let $X_1, \ldots, X_n$ be a deterministic sequence. Then a corresponding sequence of independent random variables $Y_1, \ldots, Y_n$ is generated according to a c.d. $f(y|x)$, that is, given $X_l = x$, the response $Y_l$ is distributed according to the density $f(y|x)$ which should be estimated.

In what follows, it is always assumed that predictors take values from the unit interval $[0, 1]$, which is also the support of the random predictor $X$. Further, in a fixed design case, it is assumed that $(X_1, \ldots, X_n)$ is a permutation of $(X_{(1)}, \ldots, X_{(n)})$ generated by the algorithm $X_{(0)} = 0$, $\int_{X_{(l)}}^{X_{(l+1)}} p(x)\, dx = (n+1)^{-1}$, $l = 0, 1, \ldots, n$, $X_{(n+1)} = 1$, where $p(x)$ is a positive probability density supported on $[0, 1]$. Hereafter, $p(x)$ will be referred to as the design density, regardless of an underlying design.

The considered statistical problem is to estimate the c.d. $f(y|x)$ as a bivariate function under the Mean Integrated Squared Error (MISE) criterion with the two types of loss functions defined in the last paragraph of the Introduction. Because an underlying design (fixed or random) is unknown to the statistician, a suggested estimator should be universal (not dependent on an underlying design).

We are now in position to discuss possible assumptions about the c.d. and the design density. It is traditional in the c.d. estimation literature to consider the problem as a particular example of estimation of a bivariate density, and this explains a typical assumption that an estimated c.d. $f(y|x)$ is isotropic, meaning that it is as smooth in $y$ as in $x$. Let us recall that the most popular assumption is the twofold partial differentiability of $f(y|x)$ with respect to $y$ and $x$; see the literature mentioned in the Introduction. In general, such an assumption may be reasonable for the joint density of two abstract random variables, but in a regression setting, there are obvious differences between the predictor and the response. Only as an example, which makes this point crystal clear, let us consider an additive regression model $Y = m(X) + \varepsilon$, where $\varepsilon$ is an independent error with density $q(z)$. Then $f(y|x) = q(y - m(x))$ and it is easy to realize that the smoothness of the c.d. in $y$ is dependent solely on the smoothness of $q(z)$, while the smoothness of the c.d. in $x$ depends on the smoothness of $q(z)$ and the smoothness of the underlying regression function $m(x)$. Thus, it is prudent to assume that the c.d. may be an anisotropic bivariate function whose smoothness depends on the direction; corresponding classes of such functions will be introduced in Sections 3 and 4.

**3. Sharp local-minimax lower bound for $L_2([0, 1]^2)$ loss.** The main aim of this section is to understand how an estimated c.d. together with an underlying design density affect the MISE. To explain the employed local-minimax approach (which originated in [18]), let us recall, following that article, a classical lower local-minimax bound for estimation of a univariate density $f(y)$ over the unit interval $[0, 1]$. It is assumed that the density is close to a given pivotal density $f_0(y)$. Suppose that the pivotal density is continuous and bounded below from zero on the interval $[0, 1]$ and no assumption about $f_0(y)$ for $y$ beyond the unit interval is made. Introduce a class of densities $S(m, Q, f_0, \rho) := \{f(y) : \int_{-\infty}^{\infty} f(y) \, dy = 1, (y) \geq 0, f(y) = f_0(y) + g(y), y \in [0, 1], g \in \mathcal{S}(m, Q), \sup_{y \in [0,1]} |g(y)| < \rho\}$, where $\mathcal{S}(m, Q)$ is a Sobolev class of functions $g(y)$ that are $m$-times differentiable on $[0, 1]$ and $\mathcal{S}(m, Q) := \{g(y) : g(y) = \sum_{j=1}^{\infty} \theta_j \varphi_j(y), y \in [0, 1], \sum_{j=1}^{\infty} (\pi j)^{2m} \theta_j^2 \leq Q\}$.

Following [5, 7, 12], consider estimation of a univariate density $f(y)$ based on a sample $Y_1, \ldots, Y_n$ generated according to this density, set $a_j := (\pi j)^{2m}$, define $d_Y := d_Y(f) := \int_0^1 f(y) \, dy$ as the coefficient of difficulty, introduce a positive sequence $\mu_{Yn}$ such that $d_Y n^{-1} \sum_{j=0}^{\infty} ([a_j/\mu_{Yn}]^{1/2} - a_j)_+ := Q$, and then set

$$M_n^* := d_Y n^{-1} \sum_{j=0}^{\infty} (1 - [a_j \mu_{Yn}]^{1/2})_+.$$

Pinsker [33] evaluated $M_n^*$ and showed that $M_n^* = M_n(f)(1 + o(1))$, where

(3.1) $$M_n(f) = [P(m) Q^{1/(2m+1)}][d_Y(f)/n]^{2m/(2m+1)}$$

and

(3.2) $$P(m) = (2m+1)^{1/(2m+1)}[m/(\pi(m+1))]^{2m/(2m+1)}.$$

Note that $P(m)$ and/or $P(m) Q^{1/(2m+1)}$ may be referred to as the *Pinsker constant*. Then, it is established in [18] that for a slowly vanishing positive sequence $\rho_n$, the following local-minimax lower bound holds:

(3.3) $$\inf_{\check{f}} \sup_{f \in S(m, Q, f_0, \rho_n)} M_n^{-1}(f) E_f \left\{ \int_0^1 (\check{f}(x) - f(x))^2 \, dx \right\} \geq 1 + o(1),$$

where the infimum is taken over all possible estimators $\check{f}$ based on $n$ realizations $Y_1, \ldots, Y_n$, the pivotal density $f_0(y)$ and parameters $m$, $Q$ and $\rho_n$. Moreover, this lower bound is sharp because it is attained by data-driven estimators; see [3, 5, 9].

This is the approach that we would like to take for the c.d. problem, and this is the result to match. To this end, we introduce a similar setting for conditional density estimation. Let $m_X$ and $m_Y$ be positive integers. Consider a bivariate function $g(y, x)$, $(y, x) \in [0, 1]^2$ which is $m_Y$-times differentiable with respect to $y$ and $m_X$-times differentiable with respect to $x$ (here and in what follows, partial differentiation is meant) and which belongs to a corresponding anisotropic Sobolev class,

(3.4) $$\mathcal{S}(m_Y, m_X, Q) := \left\{ g(y, x) : g(y, x) = \sum_{j,r=0}^{\infty} \theta_{jr} \varphi_j(y) \varphi_r(x), (y, x) \in [0, 1]^2, \right.$$
$$\left. \sum_{j,r=0}^{\infty} [(\pi j)^{2m_Y} + (\pi r)^{2m_X}] \theta_{jr}^2 \leq Q \right\}.$$

This Sobolev class is well known in the statistical literature; see the discussion in [25]. Let $f_0(y|x)$, $(y,x) \in (-\infty, \infty) \times [0,1]$, be a pivotal conditional density which is continuous and bounded below from zero on $[0,1]^2$, no assumption about $f_0(y|x)$ for $(y,x)$ beyond the unit square being made. Introduce a class of conditional densities $\mathcal{S}(m_Y, m_X, Q, f_0(y|x), \rho) := \{f(y|x) \colon \int_{-\infty}^{\infty} f(y|x)\,dy = 1, f(y|x) \geq 0, (y,x) \in (-\infty, \infty) \times [0,1]; f(y|u) = f_0(y|x) + g(y,x), (y,x) \in [0,1]^2; g(y,x) \in \mathcal{S}(m_Y, m_X, Q);\ \sup_{(y,x) \in [0,1]^2} |g(y,x)| < \rho\}$. Also, let $p(x)$, $\int_0^1 p(x)\,dx = 1$, be the design density which is continuous and bounded below from zero on $[0,1]$. Then, similarly to the univariate density setting, the problem is to explore a local-minimax estimation over this c.d. class. Set $a_{jr} := (\pi j)^{2m_Y} + (\pi r)^{2m_X}$, introduce the coefficient of difficulty of estimation of the c.d. over the unit square,

$$(3.5) \qquad d := d(f,p) := \int_{[0,1]^2} f(y|x) p^{-1}(x)\,dy\,dx,$$

define a positive $\eta_n$ such that

$$(3.6) \qquad dn^{-1} \sum_{j,r=0}^{\infty} ([a_{jr}/\eta_n]^{1/2} - a_{jr})_+ := Q$$

and then set

$$(3.7) \qquad R_n^*(\mathcal{S}) := dn^{-1} \sum_{j,r=0}^{\infty} (1 - [a_{jr}\eta_n]^{1/2})_+.$$

It will be shown in Section 10 that $R_n^*(\mathcal{S}) = R_n(f,p,\mathcal{S})(1 + o(1))$, where

$$(3.8) \qquad R_n(f,p,\mathcal{S}) = [P(\alpha,\beta) Q^{1/(2\tau+1)}][d(f,p) n^{-1}]^{2\tau/(2\tau+1)},$$

$\alpha = m_Y$, $\beta = m_X$, $1/(2\tau) := 1/(2\alpha) + 1/(2\beta)$, the new Pinsker constant for estimation of the c.d. on $[0,1]^2$ is

$$(3.9) \qquad P(\alpha, \beta) := \pi^{-4\tau/(2\tau+1)} [J_1(\alpha,\beta)]^{-1/(2\tau+1)} J_2(\alpha,\beta)$$

and

$$(3.10)\ \ J_1(\alpha,\beta) := \int_{\{(u,v):u^{2\alpha}+v^{2\beta}\leq 1; u,v\geq 0\}} ([u^{2\alpha}+v^{2\beta}]^{1/2} - [u^{2\alpha}+v^{2\beta}])\,dv\,du,$$

$$(3.11)\ \ J_2(\alpha,\beta) := \int_{\{(u,v)\colon u^{2\alpha}+v^{2\beta}\leq 1; u,v\geq 0\}} (1 - [u^{2\alpha}+v^{2\beta}]^{1/2})\,dv\,du.$$

We can now present a local-minimax lower bound for c.d. estimation. In what follows, $E_{(f(y|x), p(x))}\{\cdot\}$ denotes the expectation given the c.d. $f(y|x)$ and the design density $p(x)$; note that this expectation is well defined for both random and fixed design settings and we may omit the subscript whenever no confusion arises.

THEOREM 3.1.  *Random and fixed designs are considered simultaneously. Consider $L_2([0,1]^2)$ loss. Suppose that a known design density $p(x)$ is continuous and bounded below from zero on its support $[0,1]$. Then, for a slowly vanishing positive sequence $\rho_n$, the local-minimax MISE of estimation of a conditional density $f(y|x)$ satisfies the lower bound*

$$\inf_{\check{f}} \sup_{f(y|x) \in \mathcal{S}(m_Y, m_X, Q, f_0, \rho_n)} E_{(f(y|x),p(x))} \Big\{ R_n^{-1}(f, p, \mathcal{S})$$

(3.12)
$$\times \int_{[0,1]^2} (\check{f}(y|x) - f(y|x))^2 \, dy \, dx \Big\}$$

$$\geq 1 + o(1),$$

*where $R_n(f, p, \mathcal{S})$ is defined in (3.8) and the infimum is taken over all possible c.d. estimators $\check{f}(y|x)$ based on $n$ independent pairs $(Y_1, X_1), \ldots, (Y_n, X_n)$ of observations, generated according to $(f(y|x), p(x))$, as well as on the pivotal conditional density $f_0(y|x)$, the design density $p(x)$ and the parameters $m_Y$, $m_X$, $Q$ and $\rho_n$.*

Sobolev function classes are classical in the regression literature. The density estimation literature also considers smoother function classes such as analytic ones; see discussion in [5, 20, 24, 28, 36]. Thus we shall complement the class of differentiable c.d.s considered above by two classical classes of smoother functions. Let $\gamma$, $\gamma_1$ and $\gamma_2$ be positive real numbers and recall that $Q$ is a positive real number. We begin by introducing an analytic-Sobolev class of bivariate functions,

$$\mathcal{AS}(\gamma, m_X, Q) := \Big\{ g(y,x) : g(y,x) = \sum_{j,r=0}^{\infty} \theta_{jr} \varphi_j(y) \varphi_r(x), (y,x) \in [0,1]^2,$$

(3.13)
$$\theta_{jr} = \int_{[0,1]^2} g(y,x) \varphi_j(y) \varphi_r(x) \, dy \, dx,$$

$$\sum_{j,r=0}^{\infty} [(e^{\pi \gamma j} + (\pi r)^{2m_X}) I(j+r > 0)] \theta_{jr}^2 \leq Q \Big\}.$$

This class includes bivariate functions $g(y,x)$ which are analytic in $y$ and $m_X$-fold differentiable in $x$. It is also possible that the conditional density is analytic in both $y$ and $x$. Let us then define an (anisotropic) analytic class

$$\mathcal{A}(\gamma_1, \gamma_2, Q) := \Big\{ g(y,x) : g(y,x) = \sum_{j,r=0}^{\infty} \theta_{jr} \varphi_j(y) \varphi_r(x), (y,x) \in [0,1]^2,$$

(3.14)
$$\theta_{jr} = \int_{[0,1]^2} g(y,x) \varphi_j(y) \varphi_r(x) \, dy \, dx,$$

$$\sum_{j,r=0}^{\infty}[(e^{\pi\gamma_1 j}+e^{\gamma_2 r})I(j+r>0)]\theta_{jr}^2 \leq Q\bigg\}.$$

Then, analogously to the local class $\mathcal{S}(m_Y, m_X, Q, f_0(y|x), \rho_n)$ defined above, we can introduce local classes $\mathcal{AS}(\gamma, m_X, Q, f_0(y|x), \rho_n)$ and $\mathcal{A}(\gamma_1, \gamma_2, Q, f_0(y|x), \rho_n)$. For these two classes relations (3.6)–(3.7), with corresponding $a_{jr} = (e^{\pi\gamma j} + (\pi r)^{2m_X})I(j+r>0)$ and $a_{jr} = (e^{\pi\gamma_1 j} + e^{\gamma_2 r})I(j+r>0)$, imply

$$\begin{aligned}(3.15)\quad R_n(f,p,\mathcal{AS}) &= P(m_X)Q^{1/(2m_X+1)}(d(f,p)/n)^{2m_X/(2m_X+1)}\\ &\quad \times [2m_X\ln(n)/((2m_X+1)\pi\gamma)]^{2m_X/(2m_X+1)},\end{aligned}$$

where $P(m_X)$ is defined in (3.2), and

$$(3.16)\quad R_n(f,p,\mathcal{A}) = (\pi\gamma_1\gamma_2)^{-1}d(f,p)n^{-1}\ln^2(n).$$

To shorten the presentation of lower bounds, in the following proposition, we will consider these two local classes of conditional densities together.

THEOREM 3.2. *Random and fixed designs are considered simultaneously. Consider $L_2([0,1]^2)$ loss. Suppose that a known design density $p(x)$ is continuous and bounded below from zero on its support $[0,1]$. Then, for a slowly vanishing positive sequence $\rho_n$ and $\mathcal{F}$ being either $\mathcal{AS}(\gamma, m_X, Q, f_0, \rho_n)$ or $\mathcal{AS}(\gamma_1, \gamma_2, Q, f_0, \rho_n)$, the local-minimax MISE of estimation of a conditional density $f(y|x)$ satisfies the lower bound*

$$(3.17)\quad \inf_{\check{f}} \sup_{f(y|x)\in\mathcal{F}} E_{(f(y|x),p(x))}\bigg\{R_n^{-1}(f,p,\mathcal{F})\int_{[0,1]^2}(\check{f}(y|x)-f(y|x))^2\,dy\,dx\bigg\}$$
$$\geq 1 + o(1),$$

*where $R_n(f,p,\mathcal{F})$ is defined in (3.15) or (3.16) depending on the considered class $\mathcal{F}$ and the infimum is taken over all possible c.d. estimators $\check{f}(y|x)$ based on $n$ independent pairs $(Y_1, X_1), \ldots, (Y_n, X_n)$ of observations generated according to $(f(y|x), p(x))$, as well as on the pivotal conditional density $f_0(y|x)$, the design density $p(x)$ and all parameters defining the class $\mathcal{F}$.*

A plain analysis of the coefficient of difficulty $d(f,p)$ defined in (3.5) indicates that if $\int_0^1 f(y|x)\,dy \equiv 1$, then the coefficient of difficulty does not depend on the underlying c.d. $f(y|x)$. The latter is the case if $[0,1]^2$ is the support of the joint density $f(y|x)p(x)$. Let us stress that the main reason why we are considering a local-minimax is to explore how an underlying c.d. affects the coefficient of difficulty. Using a similar proof, it is straightforward

to establish that if $L_2((-\infty, \infty) \times [0, 1])$ loss is considered, then the coefficient of difficulty does not depend on an underlying c.d. since, in this case, we always have $\int_{-\infty}^{\infty} f(y|x)\,dy \equiv 1$. Due to this remark, we can omit the analysis of local-minimax MISEs for $L_2((-\infty, \infty) \times [0, 1])$ loss and consider a classical minimax approach in the next section.

**4. Sharp minimax lower bounds for $L_2((-\infty, \infty) \times [0, 1])$ loss.** The aim of this section is to find minimax lower bounds for several anisotropic classes of conditional densities and $L_2((-\infty, \infty) \times [0, 1])$ loss. The latter is of special interest in the case of unbounded responses. Because a loss function is given, no ambiguity occurs if we use identical notation for function classes defined on $[0, 1]^2$ (as in the previous section) and those defined on $(-\infty, \infty) \times [0, 1]$ (considered in this section). The motivation for this abuse of notation is that corresponding spaces are defined in such a manner that they imply the same MISE convergence for both of the considered losses and this will allow us to shorten the presentation of results.

We begin with a Sobolev anisotropic class of conditional densities,

$$(4.1) \quad \mathcal{S}(m_Y, m_X, Q) := \bigg\{ f(y|x) : f(y|x) = \sum_{r=0}^{\infty} (2\pi)^{-1} \int_{-\infty}^{\infty} h_r(u) e^{-iuy}\,du\, \varphi_r(x), \\ f(y|x) \geq 0, \int_{-\infty}^{\infty} f(y|x)\,dy \equiv 1, \\ (y, x) \in (-\infty, \infty) \times [0, 1], \\ \sum_{r=0}^{\infty} \pi^{-1} \int_0^{\infty} [u^{2m_Y} + (\pi r)^{2m_X}]|h_r(u)|^2\,du \leq Q \bigg\}.$$

To shed light on the functions $h_r(u)$, it may be helpful to note that if $h(u|x) := \int_{-\infty}^{\infty} f(y|x) e^{iyu}\,dy$ denotes the conditional characteristic function, then

$$(4.2) \quad h_r(u) := \int_0^1 h(u|x) \varphi_r(x)\,dx$$

is its $r$th Fourier coefficient and one can write $h(u|x) = \sum_{r=0}^{\infty} h_r(u) \varphi_r(x)$. The Sobolev class (4.1) contains bivariate functions $g(y, x)$, $(y, x) \in (-\infty, \infty) \times [0, 1]$, having the square integrable generalized $m_Y$-fold partial derivative with respect to $y$ and the square integrable generalized $m_X$-fold partial derivative with respect to $x$; see the discussion in [32, 36].

Another anisotropic class to consider is an analytic-Sobolev one,

$$\mathcal{AS}(\gamma, m_X, Q) := \bigg\{ f(y|x) : f(y|x) = \sum_{r=0}^{\infty} (2\pi)^{-1} \int_{-\infty}^{\infty} h_r(u) e^{-iuy}\,du\, \varphi_r(x),$$

$$f(y|x) \geq 0, \int_{-\infty}^{\infty} f(y|x)\, dy \equiv 1,$$

(4.3)

$$(y, x) \in (-\infty, \infty) \times [0, 1],$$

$$\sum_{r=0}^{\infty} \pi^{-1} \int_0^{\infty} [e^{\gamma u} + (\pi r)^{2m_X}] |h_r(u)|^2\, du \leq Q \bigg\},$$

which is an analog of (3.13). Note that this class includes, among others, classical normal, Student and Cauchy conditional densities, as well as their mixtures and one-to-one transformations which are typical in additive regression; see the discussion in [20, 36].

Finally, similarly to (3.14), we define an anisotropic analytic class

$$\mathcal{A}(\gamma_1, \gamma_2, Q) := \bigg\{ f(y|x) : f(y|x) = \sum_{r=0}^{\infty} (2\pi)^{-1} \int_{-\infty}^{\infty} h_r(u) e^{-iuy}\, du\, \varphi_r(x),$$

(4.4)
$$f(y|x) \geq 0, \int_{-\infty}^{\infty} f(y|x)\, dy \equiv 1,$$

$$(y, x) \in (-\infty, \infty) \times [0, 1],$$

$$\sum_{r=0}^{\infty} \pi^{-1} \int_0^{\infty} [e^{\gamma_1 u} + e^{\gamma_2 r}] |h_r(u)|^2\, du \leq Q \bigg\}.$$

We can now present lower minimax bounds.

THEOREM 4.1. *Random and fixed designs are considered simultaneously. Consider the case of $L_2((-\infty, \infty) \times [0, 1])$ loss. Suppose that a known design density $p(x)$ is continuous and bounded below from zero on its support $[0, 1]$. Then*

(4.5)
$$\inf_{\check{f}} \sup_{f(y|x) \in \mathcal{F}} E_{(f(y|x), p(x))} \bigg\{ \int_0^1 \bigg[ \int_{-\infty}^{\infty} (\check{f}(y|x) - f(y|x))^2\, dy \bigg] dx \bigg\}$$

$$\geq R_n(p, \mathcal{F})(1 + o(1)),$$

*where the infimum is taken over all possible estimators $\check{f}$ based on the design density $p(x)$, the class $\mathcal{F}$ and $n$ independent pairs of observations $(Y_l, X_l)$, $l = 1, \ldots, n$, generated according to $(f(y|x), p(x))$. The asymptotic minimax risk $R_n(p, \mathcal{F})$ is defined as follows. For the considered loss, the coefficient of difficulty is simplified to $d := d(p) = \int_0^1 p^{-1}(x)\, dx$, and then,*

(a) *for an anisotropic Sobolev class $\mathcal{F} = \mathcal{S}(m_Y, m_X, Q)$, the risk $R_n(p, \mathcal{F})$ is equal to the right-hand side of (3.8), with $\alpha = m_Y$, $\beta = m_X$ and $d(f, p)$ replaced by $d(p)$;*

(b) *for an analytic-Sobolev class $\mathcal{F} = \mathcal{AS}(\gamma, m_X, Q)$, the risk $R_n(p, \mathcal{F})$ is equal to the right-hand side of (3.15), with $d(f, p)$ replaced by $d(p)$;*

(c) *for an anisotropic analytic class $\mathcal{F} = \mathcal{A}(\gamma_1, \gamma_2, Q)$, the risk $R_n(p, \mathcal{F})$ is equal to the right-hand side of (3.16), with $d(f, p)$ replaced by $d(p)$.*

We have obtained lower bounds for the MISE which allow us to introduce the notion of sharp minimax estimation of the c.d. over a class $\mathcal{F}$ of c.d.s whenever the MISE of an estimator attains a corresponding lower bound.

**5. EP conditional density estimator.** The objective of this section is to suggest a conditional density estimator which is (i) adaptive to (in general unknown) design of predictors; (ii) simultaneously sharp minimax over the aforementioned anisotropic classes of conditional densities; (iii) at least univariate-rate minimax when the c.d. is a univariate density, that is, under the classical null hypothesis "the response is independent of the predictor."

The last aim makes it reasonable to rewrite a c.d. as a sum of a univariate component and a bivariate component,

$$(5.1) \qquad f(y|x) = f(y) + \psi(y, x),$$

where $\psi(y, x)$ vanishes if the response is independent of the predictor. A pair of these components is defined differently for the two studied loss functions and definitions will be presented shortly. (Let us stress that a loss function is known to the statistician, so an estimator may be chosen accordingly.) Then a blockwise-shrinkage Efromovich–Pinsker (EP) estimator will be developed for the estimation of $f(y)$ and $\psi(y, x)$.

REMARK 5.1. The interested reader can find a comprehensive discussion of the EP estimation procedure in [5]. Here we briefly recall its main idea. Suppose that a bivariate function $g(u, v)$ is estimated on a set $A$ and suppose that there exists an orthonormal basis $\{\varphi_{js}(u, v); j, s \geq 0\}$ on $A$ such that $g(u, v) = \sum_{j,s=0}^{\infty} \kappa_{js} \varphi_{js}(u, v)$, $\kappa_{js} = \int_A g(u, v) \varphi_{js}(u, v) \, du \, dv$ for $(u, v) \in A$. Then a blockwise-shrinkage EP estimator is defined as follows. All indices $(j, s)$ are divided into nonoverlapping blocks $B_k$, $k = 1, 2, \ldots$. Also, a sequence of positive thresholds $t_k$ and a cutoff $K$ are chosen. Blocks, thresholds and the cutoff may depend on the sample size $n$. An EP estimator can then be written as

$$\tilde{g}(u, v) := \sum_{k=1}^{K} \tilde{\mu}_k \sum_{(j,s) \in B_k} \tilde{\kappa}_{js} \varphi_{js}(u, v), \qquad (u, v) \in A,$$

where $\tilde{\kappa}_{js}$ is an estimator of $\kappa_{js}$, with a method of moments estimator being a typical choice, and $\tilde{\mu}_k = \mu(B_k, t_k, n, \{\tilde{\kappa}_{js}, (j, s) \in B_k\})$ is a shrinkage coefficient. Let us stress two facts about this estimator. First, neither blocks $B_k$, nor thresholds $t_k$, nor the cutoff $K$, nor the shrinkage-coefficient function $\mu$ depends on observations; instead, they are chosen a priori. Second, adaptation to the smoothness of an underlying function $g$ is achieved solely by shrinkage coefficients $\tilde{\mu}_k$ via their dependence on observations. The underlying idea of using a shrinkage procedure is to mimic Wiener's optimal shrinkage coefficient (oracle) $\mu_k^* := \sum_{(j,s)\in B_k} \kappa_{js}^2 / [\sum_{(j,s)\in B_k} (\kappa_{js}^2 + \text{Var}(\tilde{\kappa}_{js}))]$.

Let us present assumptions and notation.

ASSUMPTION 1. An estimated conditional density $f(y|x)$ belongs to a Sobolev class $\mathcal{S}(1,1,C)$, $C < \infty$, defined either in (3.4) or (4.1) for $L_2([0,1]^2)$ loss or $L_2((-\infty,\infty) \times [0,1])$ loss, respectively.

REMARK 5.2. It is convenient to define the Sobolev class (3.4) as a class of bivariate functions and the Sobolev class (4.1) as a class of conditional densities. Nonetheless, because $f(y|x)$ is the c.d., this difference plays no role in Assumption 1.

ASSUMPTION 2. The design density $p(x)$, $x \in [0,1]$, is bounded below from zero on its support $[0,1]$ and its first derivative $p^{(1)}(x)$ exists and is bounded on $[0,1]$.

NOTATION. Whenever no ambiguity may arise, sets of integration are omitted and, for instance, double integrals are taken over $[0,1]^2$ or $(-\infty,\infty) \times [0,1]$, depending on the loss under consideration. $\mathcal{N}$ denotes the set of nonnegative integers. Two given arrays of nonnegative numbers $\{0 = b_1 < b_2 < \cdots\}$ and $\{b_1' = 0, b_2' = 1 + \lfloor \ln^{3/4}(n) \rfloor, b_{s+1}' = b_s' + \lfloor b_2'(1 + 1/\ln\ln(n))^{s-2} \rfloor, s = 2, 3, \ldots\}$, will be used to define blocks and two given arrays of positive numbers, $\{t_1, t_2, \ldots\}$ and $\{t_{k\tau} := 1/\ln\ln((k+3)(\tau+3)), k, \tau = 1, 2, \ldots\}$, will denote thresholds. Two different arrays of blocks are used for estimation of univariate $f(y)$ and bivariate $\psi(y,x)$ components of $f(y|x)$; recall (5.1). The former is $\{B_k, k = 1, 2, \ldots\}$, where the blocks are either consecutive sets of nonnegative integers $B_k := \{j : b_k \leq j < b_{k+1}, j \in \mathcal{N}\}$ or intervals $B_k := [b_k, b_{k+1})$ for $L_2([0,1]^2)$ or $L_2((-\infty,\infty) \times [0,1])$ loss, respectively. The latter blocks are either sets of pairs of integers $B_{k\tau} := \{(j,r) : b_k' \leq j < b_{k+1}', b_\tau' < r \leq b_{\tau+1}', j, r \in \mathcal{N}\}$ for the $L_2([0,1]^2)$ loss or sets of mixed pairs of real and integer numbers $B_{k\tau} := \{(u,r) : u \in [b_k', b_{k+1}'), b_\tau' < r \leq b_{\tau+1}', r \in \mathcal{N}\}$ for the other loss. The corresponding lengths/cardinalities of these blocks are $L_k := b_{k+1} - b_k$ and $L_{k\tau} := (b_{k+1}' - b_k')(b_{\tau+1}' - b_\tau')$. In oracle inequalities, we shall also use so-called *adjusted lengths*

$$(5.2) \quad L_{k\tau}^* := \frac{L_{k\tau}}{\sum_{(j,r) \in B_{k\tau}} [|\int_0^1 p^{-1}(x) \varphi_{2r}(x) \, dx| + \int_0^1 |\int_0^1 f(y|x) \varphi_j(y) \, dy|^2 \, dx]}$$

or

$$
(5.3) \quad L_{k\tau}^* := L_{k\tau} \bigg/ \bigg( \sum_{r=1}^{\infty} \int_0^{\infty} I((u,r) \in B_{k\tau}) \bigg[ \bigg| \int_0^1 p^{-1}(x) \varphi_{2r}(x) \, dx \bigg| \\ + \int_0^1 \bigg| \int_{-\infty}^{\infty} e^{iuy} f(y|x) \, dy \bigg|^2 dx \bigg] du \bigg)
$$

for $L_2([0,1]^2)$ loss or the other loss, respectively, and let $L_k^* := L_k$.

We can now define EP estimators for the two losses in turn.

*EP estimator for $L_2([0,1]^2)$ loss.* Here, the c.d. $f(y|x)$ is estimated over the unit square $[0,1]^2$ and then [recalling the decomposition (5.1)]

$$
(5.4) \quad f(y) = \sum_{j=0}^{\infty} \theta_j \varphi_j(y), \qquad \theta_j := \int_{[0,1]^2} f(y|x) \varphi_j(y) \, dy \, dx.
$$

The Fourier series (5.4) implies a familiar univariate EP estimator,

$$
(5.5) \quad \tilde{f}(y) := \sum_{k=1}^{K} \tilde{\mu}_k \sum_{j \in B_k} \hat{\theta}_j \varphi_j(y),
$$

where the $\hat{\theta}_j$ are empirical Fourier coefficients,

$$
(5.6) \quad \hat{\theta}_j := n^{-1} \sum_{l=1}^{n} I(Y_l \in [0,1]) \varphi_j(Y_l) \hat{p}^{-1}(X_l),
$$

and the $\tilde{\mu}_k$ are plugged-in Wiener shrinkage coefficients,

$$
(5.7) \quad \tilde{\mu}_k := \frac{\tilde{\Theta}_k}{\tilde{\Theta}_k + \hat{d}n^{-1}} I(\tilde{\Theta}_k > t_k \hat{d} n^{-1}),
$$

where

$$
(5.8) \quad \tilde{\Theta}_k := L_k^{-1} \sum_{j \in B_k} |\hat{\theta}_j|^2 - \hat{d} n^{-1}
$$

and

$$
(5.9) \quad \hat{d} := n^{-1} \sum_{l=1}^{n} I(Y_l \in [0,1]) \hat{p}^{-2}(X_l)
$$

estimates the coefficient of difficulty $d = \int_{[0,1]^2} f(y|x) p^{-1}(x) \, dy \, dx$. The truncated-from-zero design density estimator is

$$
(5.10) \quad \hat{p}(x) := \max(1/\ln\ln(n), \tilde{p}(x)),
$$

the pivotal design density estimator $\tilde{p}(x)$ is an orthogonal series estimator

(5.11) $$\tilde{p}(x) := 1 + n^{-1} \sum_{r=1}^{n^{1/3}} \sum_{l=1}^{n} \varphi_r(X_l),$$

and the cutoff $K$ used in (5.5) is the minimal integer such that $b_{K+1} > n^{1/3} \ln \ln(n)$.

The underlying idea of EP estimation of the bivariate function $\psi(y, x)$ is based on the expansion

(5.12) $$\psi(y, x) = \sum_{j=0}^{\infty} \sum_{r=1}^{\infty} \theta_{jr} \varphi_j(y) \varphi_r(x),$$

$$\theta_{jr} := \int_{[0,1]^2} f(y|x) \varphi_j(y) \varphi_r(x)\, dy\, dx.$$

Note that the bivariate function (5.12) vanishes (as it should) if $f(y|x)$ does not depend on $x$. The corresponding bivariate blockwise-shrinkage EP estimator is then

(5.13) $$\tilde{\psi}(y, x) := \sum_{k,\tau=1}^{T} \tilde{\mu}_{k\tau} \sum_{(j,r) \in B_{k\tau}} \hat{\theta}_{jr} \varphi_j(y) \varphi_r(x),$$

(5.14) $$\hat{\theta}_{jr} := n^{-1} \sum_{l=1}^{n} I(Y_l \in [0,1]) \varphi_j(Y_l) \varphi_r(X_l) \hat{p}^{-1}(X_l),$$

(5.15) $$\tilde{\mu}_{k\tau} := \frac{\tilde{\Theta}_{k\tau}}{\tilde{\Theta}_{k\tau} + \hat{d}n^{-1}} I(\tilde{\Theta}_{k\tau} > t_{k\tau} \hat{d} n^{-1}),$$

(5.16) $$\tilde{\Theta}_{k\tau} := L_{k\tau}^{-1} \sum_{(j,r) \in B_{k\tau}} \hat{\theta}_{jr}^2 - \hat{d}n^{-1},$$

where $\hat{d}$ is defined in (5.9), $\hat{p}(x)$ in (5.10) and $T$ is the minimal integer such that $b'_{T+1} > n^{1/4} \ln \ln(n)$. The EP estimator is then defined as $\tilde{f}(y|x) := \tilde{f}(y) + \tilde{\psi}(y, x)$.

*EP estimator for $L_2((-\infty, \infty) \times [0, 1])$ loss.* Here, in representation (5.1), we have

(5.17) $$f(y) = (2\pi)^{-1} \int_{-\infty}^{\infty} h_0(u) e^{-iuy}\, dy$$

and

(5.18) $$\psi(y, x) = \sum_{r=1}^{\infty} (2\pi)^{-1} \int_{-\infty}^{\infty} h_r(u) e^{-iuy}\, du\, \varphi_r(x),$$

where

$$(5.19) \qquad h_r(u) := \int_0^1 h(u|x)\varphi_r(x)\,dx, \qquad r = 0, 1, \ldots,$$

and

$$(5.20) \qquad h(u|x) := \int_{-\infty}^{\infty} f(y|x)e^{iuy}\,dy$$

is the conditional characteristic function of $Y$ given $X = x$. This approach implies the following EP estimation of $f(y)$:

$$(5.21) \qquad \tilde{f}(y) := \pi^{-1}\int_0^{\infty} \operatorname{Re}\{\tilde{h}_0(u)e^{-iuy}\}\,du,$$

where

$$(5.22) \qquad \tilde{h}_0(u) := \sum_{k=1}^{K} \tilde{\mu}_k \hat{h}_0(u) I(u \in B_k), \qquad u \geq 0,$$

$$(5.23) \qquad \hat{h}_r(u) := n^{-1}\sum_{l=1}^{n} e^{iuY_l}\varphi_r(X_l)\hat{p}^{-1}(X_l), \qquad r = 0, 1, \ldots,$$

$$(5.24) \qquad \tilde{\mu}_k := \frac{\tilde{\Theta}_k}{\tilde{\Theta}_k + \tilde{d}n^{-1}} I(\tilde{\Theta}_k > t_k \tilde{d}n^{-1}),$$

$$(5.25) \qquad \tilde{\Theta}_k := L_k^{-1}\int_{B_k} |\hat{h}_0(u)|^2\,du - \tilde{d}n^{-1},$$

and because for the considered loss the coefficient of difficulty simplifies to $d = \int_0^1 p^{-1}(x)\,dx$, we use the estimate $\hat{p}(x)$ defined in (5.10) and then set

$$(5.26) \qquad \tilde{d} := \int_0^1 \hat{p}^{-1}(x)\,dx.$$

Further, the EP estimator of $\psi(y, x)$ is defined as

$$(5.27) \qquad \tilde{\psi}(y,x) := \pi^{-1} \sum_{k,\tau=1}^{T} \tilde{\mu}_{k\tau} \sum_{r=1}^{\infty} \int_0^{\infty} I((u,r) \in B_{k\tau})\operatorname{Re}\{\hat{h}_r(u)e^{-iuy}\}\,du\,\varphi_r(x),$$

where

$$(5.28) \qquad \tilde{\mu}_{k\tau} := \frac{\tilde{\Theta}_{k\tau}}{\tilde{\Theta}_{k\tau} + \tilde{d}n^{-1}} I(\tilde{\Theta}_{k\tau} > t_{k\tau}\tilde{d}n^{-1})$$

and

$$(5.29) \qquad \tilde{\Theta}_{k\tau} := L_{k\tau}^{-1}\sum_{r=1}^{\infty}\int_0^{\infty} I((u,r) \in B_{k\tau})|\hat{h}_r(u)|^2\,du - \tilde{d}n^{-1}.$$

The EP c.d. estimator is then defined as

$$\tilde{f}(y|x) := \tilde{f}(y) + \tilde{\psi}(y, x). \tag{5.30}$$

Note that definition (5.30) of the EP c.d. estimator is the same for both losses, but the two additive components are different. This abuse of notation will allow us to consider the two losses simultaneously.

**6. Oracle inequality.** The aim of this section is to show that the MISE of the EP estimator matches the MISE of an oracle that knows an estimated c.d. and has excellent statistical properties. As in the previous section, we are simultaneously considering two possible designs of predictors and two losses, $L_2([0,1]^2)$ and $L_2((-\infty, \infty) \times [0,1])$.

Let us introduce a blockwise-shrinkage oracle $\tilde{f}^*(y|x)$, motivated by Wiener's filter, which serves as a benchmark for the EP c.d. estimator $\tilde{f}(y|x)$. It is defined as the estimator (5.30) with estimated shrinkage coefficients replaced by coefficients depending on $(f, p)$,

$$\mu_k := \Theta_k/[\Theta_k + d(f,p)n^{-1}] \quad \text{and} \quad \mu_{k\tau} := \Theta_{k\tau}/[\Theta_{k\tau} + d(f,p)n^{-1}], \tag{6.1}$$

in place of the corresponding statistics $\tilde{\mu}_k$ defined in (5.7) or (5.24) (for the two losses) and $\tilde{\mu}_{k\tau}$ defined in (5.15) or (5.28) (for the two losses), respectively. Here and in what follows, $d(f,p)$ is defined in (3.5) for $L_2([0,1]^2)$ loss, $d(f,p) = d(p) = \int_0^1 p^{-1}(x)\,dx$ for $L_2((-\infty, \infty) \times [0,1])$ loss, and $\Theta_k$ and $\Theta_{k\tau}$ are Sobolev functionals defined as

$$\Theta_k := L_k^{-1} \sum_{j \in B_k} \theta_j^2, \qquad \Theta_{k\tau} := L_{k\tau}^{-1} \sum_{(j,r) \in B_{k\tau}} \theta_{jr}^2 \tag{6.2}$$

for $L_2([0,1]^2)$ loss and

$$\Theta_k := L_k^{-1} \int_{B_k} |h_0(u)|^2\,du,$$
$$\Theta_{k\tau} := L_{k\tau}^{-1} \sum_{r=1}^{\infty} \int_0^{\infty} I((u,r) \in B_{k\tau})|h_r(u)|^2\,du \tag{6.3}$$

for the other loss.

THEOREM 6.1. *The cases of bounded and unbounded responses, as well as the cases of the two studied losses, are considered simultaneously. Suppose that Assumptions 1 and 2 hold. Then the following oracle inequality holds for the EP estimator $\tilde{f}(y|x)$:*

$$E \int (\tilde{f}(y|x) - f(y|x))^2\,dy\,dx \leq E \int (\tilde{f}^*(y|x) - f(y|x))^2\,dy\,dx + \delta_n, \tag{6.4}$$

*where*

$$\delta_n \leq Cn^{-1}\left[\sum_{k=1}^{K}[L_k\mu_k(t_k^{1/2} + L_k^{-1/2}t_k^{-3/2}) + L_k^{-2}t_k^{-5}]\right.$$
(6.5)
$$\left. + \sum_{k,\tau=1}^{T}[L_{k\tau}\mu_{k\tau}(t_{k\tau}^{1/2} + (L_{k\tau}^*)^{-1/2}t_{k\tau}^{-3/2}) + (L_{k\tau}^*)^{-2}t_{k\tau}^{-5}]\right],$$

*and the oracle's MISE satisfies*

$$E\int(\tilde{f}^*(y|x) - f(y|x))^2\,dy\,dx$$

(6.6)
$$= n^{-1}c^*d(f,p)\left[\sum_{k=1}^{K}L_k\mu_k + \sum_{k,\tau=1}^{T}L_{k\tau}\mu_{k\tau}\right]$$

$$+ c^*\left[\sum_{k>K}L_k\Theta_k + \sum_{k,\tau=1}^{\infty}I((k,\tau)\notin[1,T]^2)L_{k\tau}\Theta_{k\tau}\right] + \delta_n^*,$$

*where $c^* = 1$ for $L_2([0,1]^2)$ loss and $c^* = \pi^{-1}$ for $L_2((-\infty,\infty)\times[0,1])$ loss and where, for any two arrays $\{\nu_k \in (0,1), k=1,2,\ldots\}$ and $\{\nu_{k\tau} \in (0,1), k,\tau = 1,2,\ldots\}$,*

$$|\delta_n^*| \leq c^*d(f,p)n^{-1}\sum_{k=1}^{K}L_k\mu_k[\nu_k + C\nu_k^{-1}\mu_k(L_k^{-1} + n^{-1/4})]$$
(6.7)
$$+ c^*d(f,p)n^{-1}\sum_{k,\tau=1}^{T}L_{k\tau}\mu_{k\tau}[\nu_{k\tau} + C\nu_{k\tau}^{-1}\mu_{k\tau}(L_{k\tau}^*)^{-1}].$$

The oracle inequality shows how well the EP estimator matches the oracle's risk. Note that it is a pointwise inequality (it is valid for a particular underlying c.d.) and it is exact (not asymptotic). The oracle inequality also allows us to establish minimax properties of the EP estimator, and this will be done in the next section.

**7. Minimax properties of the EP c.d. estimator.** The oracle inequality of Theorem 6.1 allows us to establish a number of useful minimax results. We need an extra assumption on blocks and thresholds which is common in the literature; see the discussion in [3, 5, 6].

ASSUMPTION 3. *Assume that $t_k \to 0$ and $L_{k+1}/L_k \to 1$ as $k \to \infty$ and that $\sum_{k=1}^{\infty}L_k^{-2}t_k^{-5} < \infty$.*

REMARK 7.1. To simultaneously consider the minimax approaches of Sections 3 and 4, it is assumed that for the local-minimax approach of Section 3, only an unknown additive component $g(y,x)$ of the c.d. $f(y|x) = f_0(y|x) + g(y,x)$ is estimated by the bivariate EP estimator based on empirical Fourier coefficients (5.6) and (5.14) minus corresponding Fourier coefficients of the known pivotal c.d. $f_0(y|x)$. Note that a pivotal c.d. traditionally studied in upper bounds is $f_0(y|x) = c \leq 1$, $(y,x) \in [0,1]^2$, for which there is no difference between the local-minimax and minimax EP estimators, apart from estimation of the single Fourier coefficient $\theta_0$.

THEOREM 7.1. *The cases of fixed and random designs are considered simultaneously. Let Assumptions 1–3 hold. Then for each loss, a corresponding EP c.d. estimator, defined in Section 5, is simultaneously sharp minimax over Sobolev, analytic-Sobolev and analytic classes of conditional densities considered in Sections 3 and 4, that is, the MISE of the EP c.d. estimator attains the lower bounds of Sections 3 and 4.*

We are now in position to show how well the bivariate EP estimator will perform in the case of the classical hypothesis, "the response is independent of the predictor." Under this hypothesis, suppose that an oracle $\tilde f^*(y)$ knows that the response is independent and then estimates the univariate density $f(y) = f(y|x)$ based on $n$ i.i.d. responses $Y_1, Y_2, \ldots, Y_n$. Obviously, this univariate oracle can be considered as a benchmark for any bivariate c.d. estimator given that the hypothesis is true. Our aim is to compare the bivariate EP c.d. estimator developed above with such an oracle.

We shall consider three classical classes of univariate densities. In what follows, $\alpha$ is a positive integer and $\gamma$, $Q$ and $q$ are positive real numbers. Let us begin with a class of differentiable univariate densities. For $L_2((-\infty, \infty))$ loss (recall that now a univariate density is estimated), we introduce a familiar Sobolev class $\mathcal{S}(\alpha, Q) := \{f(y) : \int_{-\infty}^{\infty} (f^{(\alpha)}(y))^2 \, dy \leq Q\} = \{h(u) : \pi^{-1} \int_0^{\infty} |u|^{2\alpha} |h(u)|^2 \, du \leq Q\}$, where $f^{(\alpha)}$ is the $\alpha$th generalized derivative and $h(u) := \int_{-\infty}^{\infty} f(y) e^{iuy} \, dy$ is the characteristic function; see [19, 32, 36]. With some obvious abuse of notation, for the case of $L_2([0,1])$ loss, we define a similar Sobolev class $\mathcal{S}(\alpha, Q) := \{f(y) : \sum_{j=1}^{\infty} (\pi j)^{2\alpha} \theta_j^2 \leq Q; \theta_0 \geq c_0 > 0, \theta_j := \int_0^{\infty} f(y) \varphi_j(y) \, dy\}$.

Let us now consider analytic densities. For the case of $L_2((-\infty, \infty))$ loss, a class of such densities was introduced in [20]: $\mathcal{A}(\gamma, Q) := \{f(y) : \pi^{-1} \int_0^{\infty} e^{\gamma u} \times |h(u)|^2 \, du < Q; h(u) = \int_{-\infty}^{\infty} e^{iuy} f(y) \, dy\}$. An $L_2([0,1])$ counterpart of this class is $\mathcal{A}(\gamma, Q) := \{f : \sum_{j=1}^{\infty} e^{\pi \gamma j} \theta_j^2 \leq Q, \theta_0 \geq c_0 > 0; \theta_j = \int_0^1 f(y) \varphi_j(y) \, dy\}$; see [5].

Finally, a bounded spectrum class is defined as $\mathcal{B}(q) := \{f(y) : h(u) = 0, |u| > q; h(u) := \int_{-\infty}^{\infty} e^{iuy} f(y) \, dy\}$ or $\mathcal{B}(q) := \{f(y) : \theta_j = 0, j > q, \theta_0 \geq c_0 >$

$0; \theta_j := \int_0^1 f(y)\varphi_j(y)\,dy\}$ for $L_2((-\infty,\infty))$ or $L_2([0,1])$ loss, respectively. A comprehensive discussion of this class can be found in [30]. Let us recall that Ibragimov and Hasminskii [24, 28] have established that the minimax rate of convergence for this class is parametric $qn^{-1}$.

Note that in all these definitions it is assumed that $f(y)$ is the density, that is, that $f(y) \geq 0$ and $\int_{-\infty}^{\infty} f(y)\,dy = 1$.

The following univariate minimax result is well known; see [5, 8].

PROPOSITION 7.1. *Suppose that the response is independent of the predictor. We are simultaneously considering the cases of $L_2([0,1])$ and $L_2((-\infty,\infty))$ losses. There exists an oracle $\check{f}^*(y)$, based on $n$ i.i.d. observations $Y_1, Y_2, \ldots, Y_n$ of the response, which is simultaneously rate minimax for bounded spectrum densities and sharp minimax for Sobolev and analytic densities. In particular, the EP univariate density estimator of [3], with blocks and thresholds $\{(L_k, t_k)\}$ satisfying Assumption 3, may serve as such an oracle and then*

$$(7.1) \quad \sup_{f \in \mathcal{B}(q)} [d(f)]^{-1} E \int (\check{f}^*(y) - f(y))^2 \, dy \leq Cqn^{-1},$$

$$(7.2) \quad \sup_{f \in \mathcal{S}(\alpha, Q)} [d(f)]^{-2\alpha/(2\alpha+1)} E \int (\check{f}^*(y) - f(y))^2 \, dy$$
$$= P(\alpha) Q^{1/(2\alpha+1)} n^{-2\alpha/(2\alpha+1)}(1 + o(1)),$$

$$(7.3) \quad \sup_{f \in \mathcal{A}(\gamma, Q)} [d(f)]^{-1} E \int (\check{f}^*(y) - f(y))^2 \, dy = (\pi \gamma n / \ln(n))^{-1}(1 + o(1)),$$

*where $P(\alpha)$ is defined in (3.2), $d(f)$ is either $\int_0^1 f(y)\,dy$ or 1 and the integrals in (7.1)–(7.3) are taken over $[0,1]$ or $(-\infty, \infty)$ for $L_2([0,1])$ loss or $L_2((-\infty, \infty))$ loss, respectively.*

We can now formulate a minimax assertion for the independent response case.

THEOREM 7.2. *The cases of fixed and random designs as well as the cases of two studied losses are considered simultaneously. Suppose that the response is independent of the predictor, that is, $f(y|x) = f(y)$, $x \in [0,1]$, and Assumptions 1–3 hold. Then the EP c.d. estimator $\tilde{f}(y|x)$ of Section 5 is simultaneously rate minimax over bounded spectrum, analytic and Sobolev classes of univariate densities, namely*

$$(7.4) \quad \sup_{f \in \mathcal{B}(q)} [d(f,p)]^{-1} E \int (\tilde{f}(y|x) - f(y))^2 \, dy\, dx \leq Cqn^{-1},$$

$$\sup_{f \in \mathcal{S}(\alpha,Q)} [d(f,p)]^{-2\alpha/(2\alpha+1)} E \int (\tilde{f}(y|x) - f(y))^2 \, dy \, dx$$
(7.5)
$$\leq P(\alpha) Q^{1/(2\alpha+1)} (n)^{-2\alpha/(2\alpha+1)} (1+o(1)),$$

$$\sup_{f \in \mathcal{A}(\gamma,Q)} [d(f,p)]^{-1} E \int (\tilde{f}(y|x) - f(y))^2 \, dy \, dx$$
(7.6)
$$\leq (\pi \gamma n / \ln(n))^{-1} (1+o(1)),$$

where $d(f,p) = \int_0^1 f(y) \, dy \int_0^1 p^{-1}(x) \, dx$ for the $L_2([0,1]^2)$ loss and $d(f,p) = \int_0^1 p^{-1}(x) \, dx$ for the $L_2((-\infty,\infty) \times [0,1])$ loss.

This theorem implies the following sharp-minimax result.

COROLLARY 7.1. *Let the assumptions of Theorem 7.2 hold. Consider the case of the uniform design density $p(x) = 1$, $x \in [0,1]$. Then the EP c.d. estimator $\tilde{f}(y|x)$ is simultaneously sharp minimax over Sobolev and analytic univariate density classes. Further, the MISE of this bivariate estimator matches the MISE of the univariate oracle $\check{f}^*(y)$ introduced in Proposition 7.1, namely,*

(7.7)
$$E \int (\tilde{f}(y|x) - f(y))^2 \, dy \, dx$$
$$= (1 + o(1)) E \int (\check{f}^*(y) - f(y))^2 \, dy + o(1) n^{-1}.$$

To avoid any possible confusion, let us explain the integrals in (7.7). Under $L_2([0,1]^2)$ loss, the left integral is taken over $[0,1]^2$, while the right one is taken over $[0,1]$ and the additive components of the EP estimator $\tilde{f}(y|x) = \tilde{f}(y) + \tilde{\psi}(y,x)$ are defined in (5.5) and (5.13). Under $L_2((\infty,\infty) \times [0,1])$ loss, the left integral in (7.7) is taken over $(-\infty,\infty) \times [0,1]$ with $y \in (-\infty,\infty)$ and $x \in [0,1]$, while the right integral is taken over $(-\infty,\infty)$. Also, in this case, the additive components of the EP estimator are defined in (5.21) and (5.27).

**8. Optimal design of predictors for c.d. estimation.** In a controlled experiment, the statistician can choose an underlying design density. Obtained results allow us to recommend a particular design density which minimizes the MISE of c.d. estimation.

It is worthwhile to begin by recalling a known result for regression function estimation. Consider a classical heteroscedastic regression $Y = m(X) + \sigma(X)\varepsilon$ with the predictor $X$ being supported on $[0,1]$ and the error $\varepsilon$ being standard normal. It is shown in [5, 13] that the MISE of regression function estimation is minimized by a design density

$$(8.1) \qquad p_r^*(x) := \sigma(x) \Big/ \int_0^1 \sigma(u)\,du.$$

At the same time, according to previous sections, to minimize the MISE of estimation of the c.d., the statistician needs to minimize the coefficient of difficulty

$$d(f,p) = \int_0^1 \left[\int_A f(y|x)\,dy\right] p^{-1}(x)\,dx,$$

where $A$ is either $[0,1]$ or $(-\infty,\infty)$, depending on the loss. A simple calculation then shows that the optimal design density for c.d. estimation is

$$(8.2) \qquad p_{c.d.}^*(x) := \left[\int_A f(y|x)\,dy\right]^{1/2} \Big/ \int_0^1 \left[\int_A f(y|u)\,dy\right]^{1/2} du.$$

In general, optimal designs (8.1) and (8.2) are different, but there is one important case where the two coincide. Consider a classical homoscedastic regression [where the scale function $\sigma(x)$ is constant] and suppose that $\int_A f(y|x)\,dy = 1$, $x \in [0,1]$ [note that the latter always holds for $A = (-\infty,\infty)$]. Then the *uniform* design is simultaneously optimal for the regression and c.d. estimation problems. Furthermore, according to Corollary 7.1, if the design is uniform, then the suggested bivariate EP estimator is sharp minimax under the hypothesis that the response is independent of the predictor. We can conclude that the uniform design has a very special place in controlled regression experiments.

Of course, in general, an underlying c.d. is unknown and cannot be used in designing an optimal experiment. A sequential design of predictors may then be a feasible option; the interested reader can find a discussion of sequential designs of predictors in [11].

## 9. Discussion.

9.1. *Effect of the design density.* The obtained theoretical results show that if the design density satisfies Assumption 2 (which is a mild assumption with the main property for the design density being differentiability), then the rate of the MISE convergence is determined solely by the smoothness of the c.d. The design density may only affect the constant of the MISE convergence via the coefficient of difficulty.

9.2. *Classical example.* Consider an additive regression $Y = m(x) + \varepsilon$, where the regression error $\varepsilon$ is independent of the predictor and its density is analytic (infinitely differentiable with the normal density, Cauchy density and mixture of normal densities being the main examples). In this case, the conditional density $f(y|x)$ is also analytic in $y$. Then, if the regression function is just $\alpha$-fold differentiable and Assumption 2 holds, for a corresponding analytic-Sobolev class the minimax MISE converges with the rate $[n/\ln(n)]^{-2\alpha/(2\alpha+1)}$. This result is both good and bad news. The good news is that up to a logarithmic penalty which is a minor factor for the curse of dimensionality, the bivariate c.d. can be estimated with the same MISE accuracy as the univariate regression function $m(x)$. Furthermore, the design density can be very rough (just differentiable), even in the case of an analytic error density. The bad news is that even if the error density is analytic and can then, according to [8], be estimated with the MISE of order $\ln(n)/n$ (i.e., with almost parametric accuracy), the MISE of c.d. estimation is primarily defined by the smoothness of the regression function and thus may be dramatically larger than the MISE of error density estimation.

9.3. *Fixed design.* Fixed design regression is the classical setting in applied regression analysis; see the discussion in [31]. For this setting, definition (1.1), which has been the key for the random design case, is not valid. Fortunately, this paper shows that a design affects neither lower bounds, nor upper bounds, nor the minimax data-driven EP estimation procedure, nor oracle inequalities.

9.4. *Dimension reduction.* A traditional null hypothesis in regression analysis is that a response and a predictor are unrelated, that is, $f(y|x) = f(y)$. In this case, the accuracy of estimation under an MISE criterion must be dramatically better because the estimated function is univariate and no curse of dimensionality occurs. It is established that the EP estimator provides an optimal univariate accuracy of estimation when the null hypothesis is valid and thus solves the classical dimension reduction problem.

9.5. *C.d. estimation in other settings.* The "regression" methodology thus far developed can be used in other classical settings, for instance in the popular time series one; see the discussion about this setting in [15]. The main complication here is that pairs of observations are no longer independent; at the same time, the setting is simpler because covariates cannot be deterministic. There are also many interesting expansions in the regression setting considered. For instance, the predictor can be a vector and the covariates may be qualitative and quantitative; see the discussion of such a setting in [21]. Some new results for this setting can be found in [10].

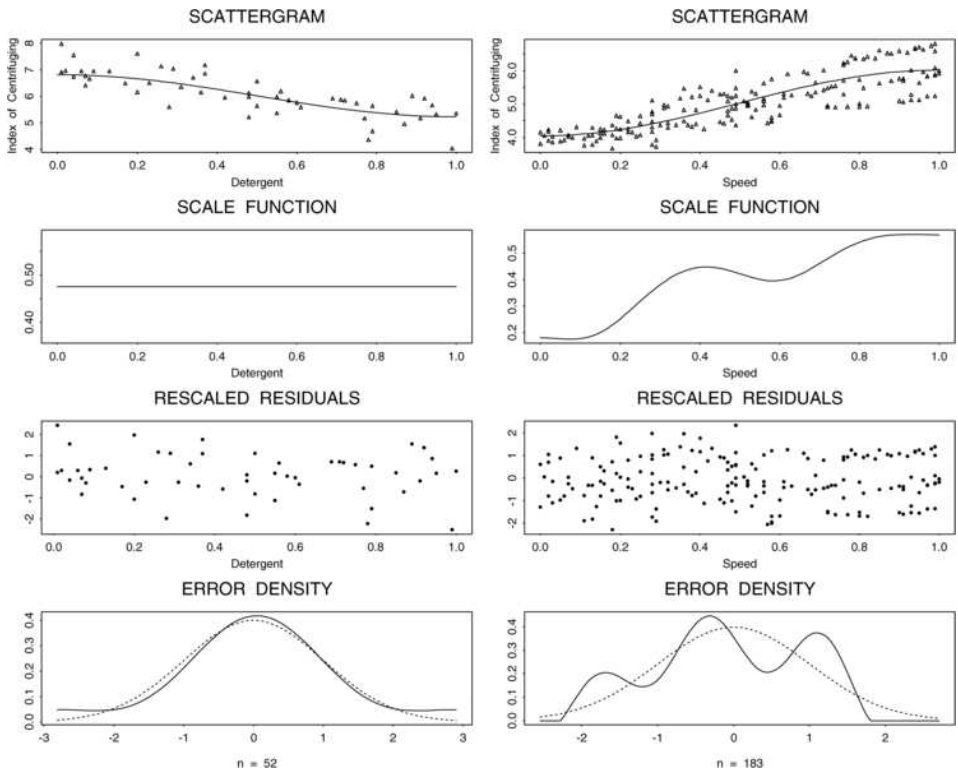

Fig. 1. *Standard nonparametric regression analysis of two real datasets. The top-left diagram exhibits a scattergram of 52 observations showing a relationship between the amount of a detergent in a sludge and an index of centrifuging. The top-right diagram shows a scattergram of 183 observations showing a relationship between speed of rotation and an index of centrifuging. Scattergrams are overlaid by nonparametric regression estimates. Dotted lines in the bottom diagrams show the standard normal density.*

9.6. *Minimax paradigm.* This paper uses a classical minimax approach: an estimator must be minimax whenever $Y$ and $X$ are dependent (the c.d. is a bivariate function) and then, if $Y$ and $X$ are independent (the c.d. is a univariate function), it is desirable that the estimator be also minimax over univariate estimators/oracles. Note that the priorities are reversed in the dimension reduction literature.

9.7. *Small datasets.* Let us begin by exploring two datasets collected by BIFAR, a company with interests in waste water treatment; the interested reader can find a complete account of these experiments in [9]. In what follows, freely available software from [5] is used; recall that it is based on mimicking EP estimators. Two columns of diagrams in Figure 1 exhibit a standard nonparametric regression analysis of two different datasets. The top-left diagram shows a scattergram with a pronounced regression function. Diagrams below it indicate that the regression is homoscedastic with normal regression errors. This dataset is a textbook example, where the regression function allows one to quantify the impact of the predictor on the response. The dataset analyzed in the right column of Figure 1 is more complicated, due to the stepwise shape of the scale function and the multimodal marginal density of regression errors; thus, let us look at the c.d. estimate shown in Figure 2. The estimate exhibits pronounced ridges and large valleys. There are several interesting features of the exhibited ridges: they are almost parallel to the speed-axis; they rise and then collapse over the speed range; ridges with larger speeds apparently have larger indices; the number of pronounced ridges increases from 1 to 3 over the range of speed. The interested reader can now return to Figure 1 and understand why the scale and error density estimates have those interesting shapes.

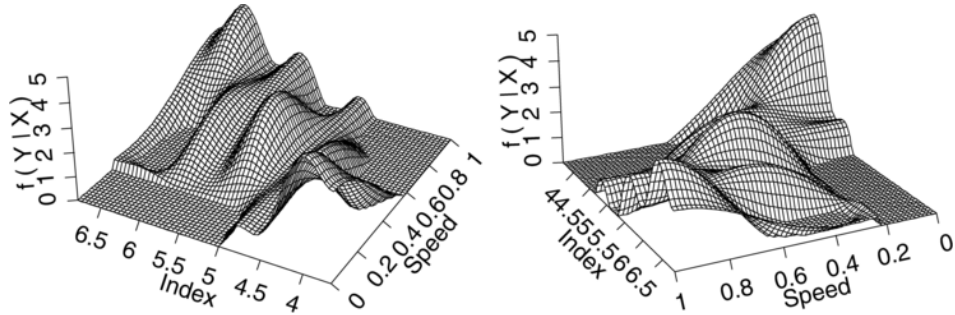

Fig. 2. *Conditional density estimate, multiplied by a factor of 3, for the speed-index dataset exhibited in the top-right diagram in Figure 1; two views are shown.*

Is it possible that estimates for the second dataset are just products of a "spurious" realization of a classical regression indicated in the first example? Let us check this by intensive Monte Carlo simulations based on the regression model for the first example and $n = 183$. Visual analysis of 500 c.d. estimates revealed that only 32 of those estimates exhibited more than one ridge, and in none of those cases did the error density estimate exhibit more than one mode. In other words, none of the Monte Carlo simulations revealed the pattern observed for the second experiment. Moreover, the author analyzed two more experiments, identical to the second one, and they exhibited similar patterns for the error density and the c.d. These results show that a "spurious" nature of the estimates is unlikely.

Let us also present results of an interesting Monte Carlo study conducted under the null hypothesis "the response and the predictor are independent." Suppose that $Y$ and $X$ are independent, $Y$ is standard normal and $X$ is standard uniform. In this study the bivariate EP c.d. estimator is compared with two univariate kernel oracles: (i) a super-oracle which knows that $Y$ and $X$ are independent and that the estimated univariate c.d. $f(y|x) \equiv f(y)$ is standard normal and which uses a Gaussian kernel with the optimal (for the underlying standard normal density) bandwidth (see [5], page 358); (ii) a sub-oracle which knows that $Y$ and $X$ are independent but does not know the density of $Y$ and which uses a Gaussian kernel, but where choice of bandwidth is done by the S-PLUS function *density*. 500 simulations were conducted for each sample size, medians of ratios of empirical ISE's of the nonparametric estimate to empirical ISE's of oracles then being calculated. For sample sizes 50, 100, 150, 200 and 300, the corresponding medians were (the numerator presenting a median ratio for the super-oracle and the denominator for the sub-oracle) 6.6/0.83, 2.21/0.31, 1.95/0.37, 2.27/0.34 and 2.62/0.52. As we see, the c.d. estimator cannot match the super-oracle which knows the underlying c.d., but it performs comparatively well when $n \geq 100$. At the same time, it outperforms the sub-oracle.

## 10. Proofs.

*Proof of Theorem* 3.1. We begin by dividing the unit square $[0,1]^2$ into $s^2$ subsquares, where the known densities $p(x)$ and $f_0(y|x)$ are approximated by constants. Lower bounds are then established for each subsquare and the total is evaluated; this is the plan of the proof. Also, whenever possible, random and fixed designs will be considered simultaneously.

Set $s = 1 + \lfloor \ln(\ln(n+20)) \rfloor$ and define $\mathcal{H}_s = \{f : f(y|x) = f_0(y|x) + [\sum_{k,r=0}^{s-1} f_{(kr)}(y|x) - \sum_{k,r=0}^{s-1} \int_0^1 f_{(kr)}(z|x)\,dz] I((y,x) \in [0,1]^2), f_{(kr)}(y|x) \in \mathcal{H}_{skr}, f(y|x) \geq 0\}$. The function classes $\mathcal{H}_{skr}$ are defined as follows. Let $\tilde{\phi}(y) := \phi(n,y)$ be a sequence of flat-top nonnegative kernels defined on the real line such that for a given $n$, the kernel is zero beyond $(0,1)$, it is $m_Y$-fold continuously differentiable on $(-\infty, \infty)$, $0 \leq \tilde{\phi}(y) \leq 1$, $\tilde{\phi}(y) = 1$ for $2(\ln(n))^{-2} \leq y \leq 1 - 2(\ln(n))^{-2}$ and its $l$th derivative satisfies $\max_y |\tilde{\phi}^{(l)}(y)| \leq C(\ln(n))^{2l}$, $l = 1, \ldots, m_Y$. For instance, such a kernel may be constructed using the so-called mollifiers, discussed in [5]. Let $\tilde{\phi}_{sk}(y) := \tilde{\phi}(sy - k)$. Analogously define $\hat{\phi}_{sr}(x)$, with $m_X$ replacing $m_Y$. Set $\varphi_{skj}(y) := s^{1/2} \varphi_j(sy - k)$. For a $(k,r)$th subsquare, $0 \leq k, r \leq s-1$, define $\phi_{skr}(y,x) := \tilde{\phi}(sy-k)\hat{\phi}(sx-r)$, $\varphi_{skrjt}(y,x) := \varphi_{skj}(y)\varphi_{srt}(x)\phi_{skr}(y,x)$, $f_{[kr]}(y|x) := \sum_{(j,t)\in T(s,k,r)} \nu_{skrjt} \varphi_{skrjt}(y,x)$ and $f_{(kr)}(y|x) := f_{[kr]}(y|x)\phi_{skr}(y,x)$. The set $T(s,k,r)$ of pairs $(j,i)$ is the difference between two sets defined as follows. Let $\eta_n(Q)$ be defined by means of the relation $\sum_{j,t \geq 0} ([a_{jt}/\eta_n(Q)]^{1/2} - a_{jt})_+ := nd^{-1}Q$ with $d$ defined in (3.5), $a_{jt} = 1 + (\pi j)^{2m_Y} + (\pi t)^{2m_X}$ and $(x)_+ = \max(0,x)$. Then the larger set is $\{(j,t) : a_{jt} \leq [\eta_n(Q_{skr})]^{-1/2}\}$ with $Q_{skr} := Q(1-1/s)(\overline{I_s^{-1}} I_{skr})^{-1}$, where $I_{skr} := p(rs^{-1})/f_0(ks^{-1}|rs^{-1})$, $\overline{I_s^{-1}} = \sum_{k,r=0}^{s-1}(1/I_{skr})$. The smaller set consists of pairs $(j,t)$ such that $\max(j,t) \leq \ln^s(n)$. We can now define $\mathcal{H}_{skr} :=$

$\{f_{(kr)}(y|x) : \sum_{(j,t) \in T(s,k,r)} [1 + (\pi sj)^{2m_Y} + (\pi st)^{2m_X}] \nu^2_{skrjt} \leq Q_{skr}, |f_{[kr]}(y|x)|^2 \leq s^4 \ln(n) R_n\}$, where $R_n := R_n(f_0, p, \mathcal{S})$ is defined in (3.8).

Let us verify that for sufficiently large $n$, we have $\mathcal{H}_s \subset \mathcal{S}(m_Y, m_X, Q, f_0, \rho_n)$. The definition of the flat-top kernel implies that $f(y|x) - f_0(y|x)$, $(y,x) \in [0,1]^2$, is $m_Y$-fold differentiable with respect to $y$ and $m_X$-fold differentiable with respect to $x$. Second, let us verify that for $f \in \mathcal{H}_s$, the difference $f(y|x) - f_0(y|x)$ belongs to $\mathcal{S}(m_Y, m_X, Q)$. Set $m = m_Y$ and begin with the differentiation with respect to $y$; in several of the following lines we use the notation $\psi^{(l)}(y,x) := \partial^l \psi(y,x)/\partial y^l$. By the Leibniz rule, $(f_{[kr]}(y|x)\phi_{skr}(y,x))^{(m)} = \sum_{l=0}^{m} \mathbf{C}_l^m f_{[kr]}^{(m-l)}(y|x) \phi_{skr}^{(l)}(y,x)$, where $\mathbf{C}_l^m := m!/((m-l)!l!)$. Note that for $0 < l \leq m$, we have $(\phi_{skr}^{(l)}(y,x))^2 \leq C(s(\ln(n))^2)^{2l}$ and for $f_{(kr)} \in \mathcal{H}_{skr}$,

$$
(10.1) \quad \begin{aligned}
&\int_{[0,1]^2} [f_{[kr]}^{(m-l)}(y|x) \phi_{skr}^{(l)}(y,x)]^2 \, dx\, dy \\
&\leq Cs^{2l} \ln^{4l}(n) \int_{k/s}^{(k+1)/s} \left( \int_{r/s}^{(r+1)/s} [f_{[kr]}^{(m-l)}(y|x)]^2 \, dx \right) dy \\
&\leq Cs^{2l} \ln^{4l}(n) \sum_{(j,t) \in T(s,k,r)} j^{2(m-l)} \nu^2_{skrjt} \\
&\leq C \ln^{4m+1}(n) \max_{(j,t) \in T(s,k,r)} \frac{j^{2(m-l)}}{1 + j^{2m} + t^{2m_X}} Q_{srk} = o(1) \ln^{-2}(n) Q_{skr}.
\end{aligned}
$$

In the last inequality we used the definition of $\mathcal{H}_{skr}$ and the assumption that $\min(j,t : (j,t) \in T_{skr}) > \ln^s(n)$. A similar conclusion can be arrived at for the derivatives with respect to $x$. Then, using Parseval's identity, we can write for $f_{(kr)} \in \mathcal{H}_{skr}$,

$$
(10.2) \quad \begin{aligned}
&\int_{[0,1]^2} [f_{[kr]}^2(y|x) + (\partial^{m_Y} f_{[kr]}(y|x)/\partial y^{m_Y})^2 \\
&\qquad + (\partial^{m_X} f_{[kr]}(y|x)/\partial x^{m_X})^2] \phi_{skr}^2(y,x)\, dx\, dy \\
&\leq \sum_{(j,t) \in T(s,k,r)} [1 + (\pi sj)^{2m_Y} + (\pi st)^{2m_X}] \nu^2_{skrjt} \leq Q_{skr}.
\end{aligned}
$$

Using this, the fact that the function $\sum_{k,r=1}^{s-1} f_{kr}(y|x)$ and its corresponding derivatives are zero at the boundary of $[0,1]^2$, Proposition 1 of [7] and the fact that $\sum_{k,r=0}^{s-1} Q_{skr} = Q(1 - s^{-1})$, we can conclude that $\sum_{k,r=1}^{s-1} f_{(kr)}(y|x) \in S(m_X, m_Y, Q(1 - s^{-1}))$. We are left with the verification that a function $g_s(x) := \sum_{k,r=1}^{s-1} \int_0^1 f_{(k,r)}(y|x)\, dy$ belongs to $S(m_X, m_Y, o(1)s^{-1})$. Write for

$f_{(kr)} \in \mathcal{H}_{skr}$,

$$g_s(x) = \sum_{k,r=0}^{s-1} \int_{k/s}^{(k+1)/s} f_{[kr]}(y|x) \phi_{skr}(y,x) \, dy$$

$$= \sum_{k,r=0}^{s-1} \int_{k/s}^{(k+1)/s} f_{[kr]}(y|x) [1 - \phi_{skr}(y,x)] \, dy,$$

where we use $\int_{k/s}^{(k+1)/s} f_{[kr]}(y|x) \, dy \equiv 0$. We then get

$$\int_0^1 [g_s^2(x) + (g_s^{(m_X)}(x))^2] \, dx$$

$$\leq o(1) \ln^{-2}(n) + \int_0^1 \left[ \sum_{k,r=0}^{s-1} \int_{k/s}^{(k+1)/s} f_{[kr]}^{(m_X)}(y|x)(1 - \phi_{skr}(y,x)) \, dy \right]^2 dx$$

$$= o(1) \ln^{-2}(n).$$

This verifies that $\mathcal{H}_s \subset \mathcal{S}(m_Y, m_X, Q, f_0, \rho_n)$ whenever $\rho_n$ vanishes slowly.

Let us now establish a lower bound for $f \in \mathcal{H}_s$ and any estimate $\hat{f}_n(y|x)$. Denote $\hat{f}(y|x) =: f_0(y|x) + \tilde{f}(y|x)$, $b_s(x) := \sum_{k,r=0}^{s-1} \int_0^1 f_{(kr)}(u|x) \, du = \sum_{k,r=0}^{s-1} \int_0^1 f_{[kr]}(x)(1 - \phi_{skr}(u,x)) \, du$ and note that for $f(y|x) \in \mathcal{H}_s$ and any $\gamma > 0$,

$$\int_{k/s}^{(k+1)/s} \int_{r/s}^{(r+1)/s} (\hat{f}(y|x) - f(y|x))^2 \, dx \, dy$$

$$= \int_{k/s}^{(k+1)/s} \int_{r/s}^{(r+1)/s} (\tilde{f}(y|x) - f_{(kr)}(y|x) + b_s(x))^2 \, dx \, dy$$

$$\geq (1-\gamma) \int_{k/s}^{(k+1)/s} \int_{r/s}^{(r+1)/s} (\tilde{f}(y|x) - f_{[kr]}(y|x))^2 \, dx \, dy$$

$$- \gamma^{-1} \int_{k/s}^{(k+1)/s} [f_{[kr]}(y|x)(1 - \phi_{skr}(y,x)) + b_s(x)]^2 \, dx$$

$$\geq (1-\gamma) \int_{k/s}^{(k+1)/s} \int_{r/s}^{(r+1)/s} (\tilde{f}(y|x) - f_{[kr]}(y,x))^2 \, dx \, dy$$

$$+ o(1) \gamma^{-1} (\ln(n))^{-1/2} R_n.$$

Then set $\gamma = s^{-1}$ and write

$$\sup_{f \in \mathcal{S}(m_Y, m_X, Q, f_0, \rho)} E \left\{ \int_{[0,1]^2} (\hat{f}(y|x) - f(y|x))^2 \, dx \, dy \right\}$$

$$\geq \sup_{f \in \mathcal{H}_s} E\left\{ \int_{[0,1]^2} (\hat{f}(y|x) - f(y|x))^2 \, dx \, dy \right\}$$

(10.3)
$$= \sup_{f \in \mathcal{H}_s} \sum_{k,r=0}^{s-1} E\left\{ \int_{k/s}^{(k+1)/s} \int_{r/s}^{(r+1)/s} (\hat{f}(y|x) - f(y|x))^2 \, dx \, dy \right\}$$

$$\geq (1 - s^{-1}) \sum_{k,r=0}^{s-1} \sup_{f \in \mathcal{H}_{skr}} \sum_{(j,t) \in T(s,k,r)} E\{(\tilde{\nu}_{skrjt} - \nu_{skrjt})^2\} + o(1) R_n$$

$$=: (1 - s^{-1}) \sum_{k=0}^{s-1} A_{kr} + o(1) R_n,$$

where $\tilde{\nu}_{skrjt} := \int_{k/s}^{(k+1)/s} \int_{r/s}^{(r+1)/s} \tilde{f}(y|x) \varphi_{skrjt}(y,x) \, dx \, dy$. As we see, the original problem is converted into the problem of finding lower bounds for terms $A_{kr}$ corresponding to a subsquare; recall that our underlying idea has been to approximate the known conditional density $f_0(y|x)$ and the univariate density $h(x)$ by constant functions on each subsquare. We continue with the following steps. First, we introduce an array of independent normal random variables $\zeta_{skrjt}$ with zero mean and variance $(1 - \gamma_n)\nu_{skrjt}^2$, where the positive sequence $\gamma_n$ tends to zero as slowly as desired. We then introduce a stochastic process $f^*(y|x)$, defined as the $f(y|x) \in \mathcal{H}_s$ previously studied, but with random $\zeta_{skrjt}$ used in place of fixed and known $\nu_{skrjt}$. The idea of considering such a stochastic process was suggested in [33], and following along the lines of the establishment of (A.18) in that article, we obtain

(10.4) $\quad P((f^*(y|x) - f_0(y|x)) \in \mathcal{S}(m_Y, m_X, Q)) = 1 + o(1).$

Now let us additionally suppose that $\nu_{skrjt}^2 \leq sn^{-1}$. It is then easily verified that

$$\sum_{(j,t) \in T(s,k,r)} \sup_{y,x} [\nu_{skrjt} \varphi_{skrjt}(y,x)]^2 \leq Cs^3 R_n.$$

Further, we can introduce a similarly defined stochastic process $f^*_{[kr]}$. This, together with Theorem 6.2.3 in [29], implies the inequality

$$P\left( \sup_{(y,x) \in [0,1]^2} |f^*_{[kr]}(y|x)|^2 \leq s^4 \ln(n) R_n \right) \geq 1 - |o(1)|s^{-2}.$$

Our next step is to compute the classical parametric Fisher information for $f \in \mathcal{H}_s$. Here, different calculations are needed for random and fixed designs. Let us begin with the former one where observations are i.i.d. pairs $(Y_l, X_l)$, $l = 1, \ldots, n$, and thus the Fisher information of $n$ pairs is $n$ times the Fisher information of a single pair. For a parameter $\nu_{skrjt}$, the "individual" Fisher information is

(10.5) $\quad I_{skrjt} := E_{(f_0, p)} \{ [\partial \ln(f(Y|X) p(X)) / \partial \nu_{skrjt}]^2 \}.$

Note that
$$\frac{\partial \ln f(y|x)}{\partial \nu_{skrjt}}$$

(10.6)
$$= \left[\partial \ln\left(\left[f_0(y|x) + \sum_{k,r=0}^{s-1} f_{(kr)}(y|x)\right.\right.\right.$$
$$\left.\left.\left. - \sum_{k,r=0}^{s} \int_0^1 f_{(kr)}(z|x)\, dz\right] I((y,x) \in [0,1]^2)\right)\right] \bigg/ \partial \nu_{skrjt}$$
$$= \frac{\varphi_{skrjt}(y,x) - \int_0^1 \varphi_{skrjt}(y,z)\, dz}{f(y|x)} I((y,x) \in [0,1]^2).$$

Recall that $f_0(y|x)p(x)$ is continuous on the unit square and write

(10.7) $\quad I_{skrjt} = \int_{[0,1]^2} f_0(y|x) p(x) \left[\frac{\varphi_{skrjt}(y,x) - \int_0^1 \varphi_{skrjt}(z,x)\, dz}{f(y|x)}\right]^2 dx\, dy.$

Further,
$$\int_{[0,1]^2} f_0(y|x) p(x) \left[\frac{\varphi_{skj}(y)\varphi_{srt}(x)\phi_{skr}(y,x)}{f(y|x)}\right]^2 dx\, dy$$
$$= \int_{k/s}^{(k+1)/s} \int_{r/s}^{(r+1)/s} f_0(y|x) p(x) \frac{\varphi_{skj}^2(y)\varphi_{srt}^2(x)}{f^2(y,x)}\, dx\, dy$$
$$+ \int_{k/s}^{(k+1)/s} \int_{r/s}^{(r+1)/s} f_0(y|x) p(x)$$
$$\times \frac{\varphi_{skj}^2(y)\varphi_{srt}^2(x)[\phi_{skr}^2(x,y) - 1]}{f^2(y|x)}\, dx\, dy$$

(10.8)
$$= \int_{k/s}^{(k+1)/s} \int_{r/s}^{(r+1)/s} [f_0(ks^{-1}|rs^{-1})p(ks^{-1}) + o(1)]$$
$$\times \frac{\varphi_{skj}^2(y)\varphi_{srt}^2(x)}{f_0^2(ks^{-1}|rs^{-1})(1+o(1))}\, dx\, dy$$
$$+ o(1)\ln^{-1}(n) = \frac{p(rs^{-1})}{f_0(ks^{-1}|rs^{-1})}(1 + o(1))$$
$$= I_{skr}(1 + o(1)).$$

Note that here $o(1) \to 0$ as $n \to \infty$ uniformly over the considered $(k,j,t)$. Also, for $j > 0$ and all sufficiently large $n$, we obtain
$$\int_{[0,1]^2} f_0(y|x) p(x) \left[\frac{\int_0^1 \varphi_{skj}(z)\varphi_{srt}(x)\phi_{skr}(z,x)\, dz}{f(y|x)}\right]^2 dx\, dy$$

$$
\begin{aligned}
&\leq C \int_0^1 \left[ \int_0^1 \varphi_{skj}(z) \varphi_{srt}(x) \phi_{skr}(z,x) \, dz \right]^2 dx \\
&\leq C \int_0^1 \left[ \int_{k/s}^{(k+1)/s} \varphi_{skj}(z) \varphi_{srt}(x) \, dz \right]^2 dx \\
&\quad + Cs^2 \int_0^1 \left[ \int_0^1 (1 - \varphi_{skr}(z,x)) \, dz \right]^2 dx \leq C \ln^{-3}(n).
\end{aligned}
\tag{10.9}
$$

Combining the obtained results in (10.5), we get $I_{skrjt} = I_{skr}(1 + o(1))$ with $o(1) \to 0$ as $n \to \infty$ uniformly over the considered $(k, r, j, t)$. Now let us calculate Fisher information for the fixed design case. Here, observations are pairs $(Y_l, X_l)$, $l = 1, \ldots, n$, where the predictors are deterministic and the responses are independent but not identically distributed. Without loss of generality, we can assume that $X_1 < X_2 < \cdots < X_n$. Note that the Fisher information of $n$ pairs is equal to the sum of the "individual" Fisher information values. Let us calculate this "individual" information for a pair $(Y_l, X_l)$ with respect to the parameter $\nu_{skrjt}$,

$$
I_{skrjt}(l) := E_{(f_0, p)}[\partial \ln(f(Y|X_l))/\partial \nu_{skrjt}]^2. \tag{10.10}
$$

Use of a calculation similar to (10.6)–(10.9) shows that $I_{skrjt}(l) = f_0^{-1}(ks^{-1}|rs^{-1}) \varphi_{srt}^2(x_l)(1 + o(1))$ if $X_l \in [rs^{-1}, (r+1)s^{-1})$ and that it is zero otherwise. This yields

$$
\begin{aligned}
\sum_{l=1}^n I_{skrjt}(l) &= f_0^{-1}(ks^{-1}|rs^{-1}) s^{-1} \sum_{\{l : X_l \in [rs^{-1}, (r+1)s^{-1}), 1 \leq l \leq n\}} \varphi_{srt}^2(X_l)(1 + o(1)) \\
&= f_0^{-1}(ks^{-1}|rs^{-1}) s^{-1} \\
&\quad \times \sum_{\{l : X_l \in [rs^{-1}, (r+1)s^{-1}), 1 \leq l \leq n\}} (X_{l+1} - X_l) \varphi_{srt}^2(X_l) p(rs^{-1}) n(1 + o(1)) \\
&= n f_0^{-1}(ks^{-1}|rs^{-1}) p(rs^{-1})(1 + o(1)) = n I_{skr}(1 + o(1)).
\end{aligned}
$$

We can conclude that asymptotically the average Fisher information is the same for both designs. With this remark in mind, we can again continue our analysis of both cases simultaneously.

We are now evaluating $\eta_n$ and $R_n^*$ as defined in (3.6)–(3.7). Set $\alpha := m_Y$, $\beta := m_X$, $N := 1/\eta_n$ and rewrite (3.6) as

$$
\sum_{\{(j,t): 0 < a_{jt} < N\}} [(a_{jt} N)^{1/2} - a_{jt}] = Q d^{-1} n. \tag{10.11}
$$

The sum in (10.11) can be approximated for large $N$ (or equivalently for large $n$) by the integral

$$
G_N := \int_{\{(y,x) \,:\, (\pi y)^{2\alpha} + (\pi x)^{2\beta} \leq N; y, x > 0\}} ([(\pi y)^{2\alpha} + (\pi x)^{2\beta}]^{1/2} N^{1/2}
$$

$$- [(\pi y)^{2\alpha} + (\pi x)^{2\beta}]) \, dx \, dy.$$

Let us apply the change of variables $u = \pi y N^{-1/(2\alpha)}$ and $v = \pi x N^{-1/(2\beta)}$. Then

$$G_N = (\pi)^{-2} N^{1/(2\alpha)} N^{1/(2\beta)} N$$
$$\times \int_{\{(u,v): u^{2\alpha} + v^{2\beta} \leq 1; u,v > 0\}} ([u^{2\alpha} + v^{2\beta}]^{1/2} - [u^{2\alpha} + v^{2\beta}]) \, dv \, du.$$

This yields $\eta_n = ([Q\pi^2 J_1^{-1}(\alpha, \beta)][d^{-1}n])^{-2\tau/(2\tau+1)}(1+o(1))$. To evaluate $R_n^*$, we again approximate the sum in (3.7) by a corresponding integral and then employ the change of variables described above,

$$G'_N := \int_{\{(y,x): (\pi y)^{2\alpha} + (\pi x)^{2\beta} \leq N; y,x > 0\}} (1 - [(\pi y)^{2\alpha} + (\pi x)^{2\beta}]^{1/2} N^{-1/2}) \, dx \, dy$$
$$= (\pi)^{-2} N^{1/(2\tau)} \int_{\{(u,v): u^{2\alpha} + v^{2\beta} \leq 1; u,v > 0\}} (1 - [u^{2\alpha} + v^{2\beta}]^{1/2}) \, dv \, du.$$

This implies that $R_n^* = P(\alpha, \beta) Q^{1/(2\tau+1)} (d/n)^{2\tau/(2\tau+1)} (1 + o(1))$.

We have established all propositions to proceed along the lines of the proof of Theorem 1 in [4]. This yields that uniformly over $k, r \in \{0, 1, \ldots, s-1\}$,

$$\inf A_{kr} \geq (s^{-4\tau} Q_{skr})^{1/(2\tau+1)} (n I_{skr})^{-2\tau/(2\tau+1)} P(m_Y, m_X)(1 + o(1)),$$
(10.12)

where the infimum is over all possible nonparametric estimates of $f$ considered in the theorem. Recalling the definition of $Q_{skr}$ and the fact that $s = s(n) \to \infty$, $n \to \infty$, we get

$$\inf \sum_{k,r=0}^{s-1} A_{kr} \geq P(m_Y, m_X) Q^{1/(2\tau+1)} n^{-2\tau/(2\tau+1)} s^{-4\tau/(2\tau+1)}$$
(10.13)
$$\times \left[ \sum_{k,r=0}^{s-1} (\overline{I_s^{-1}} I_{skr})^{-1/(2\tau+1)} I_{skr}^{-2\tau/(2\tau+1)} \right] (1 + o(1)).$$

Further,

$$\sum_{k,r=0}^{s-1} (\overline{I_s^{-1}} I_{skr})^{-1/(2\tau+1)} I_{skr}^{-2\tau/(2\tau+1)} = (\overline{I_s^{-1}})^{-1/(2\tau+1)} \sum_{k,r=0}^{s-1} I_{skr}^{-1}$$
(10.14)
$$= (\overline{I_s^{-1}})^{2\tau/(2\tau+1)}$$
$$= \left[ \sum_{k,r=0}^{s-1} \frac{f_0(ks^{-1}|rs^{-1})}{p(rs^{-1})} \right]^{2\tau/(2\tau+1)}.$$

Using our assumption about continuity of $f_0(y|x)$ and $p(x)$ on the unit square, we obtain

$$s^{-4\tau/(2\tau+1)}\left[\sum_{k,r=0}^{s-1}\frac{f_0(ks^{-1}|rs^{-1})}{p(rs^{-1})}\right]^{2\tau/(2\tau+1)}$$

$$=\left[s^{-2}\sum_{k,r=0}^{s-1}\frac{f_0(ks^{-1}|rs^{-1})}{p(rs^{-1})}\right]^{2\tau/(2\tau+1)}$$

$$=\left[\int_{[0,1]^2}\frac{f_0(y|x)}{p(x)}\,dx\,dy\right]^{2\tau/(2\tau+1)}(1+o(1)).$$

We conclude that

$$\inf\sum_{k,r=0}^{s-1}A_{kr}\geq P(m_Y,m_X)Q^{1/(2\tau+1)}$$

$$\times\left[n^{-1}\int_{[0,1]^2}\frac{f_0(y|x)}{p(x)}\,dx\,dy\right]^{2\tau/(2\tau+1)}(1+o(1)).$$

This, together with (10.3), verifies Theorem 3.1.

Proofs of the lower bounds in Theorems 3.2 and 4.1 are similar. Proofs of the upper bounds can be found in the technical report [9].

**Acknowledgments.** The author is grateful for the helpful and constructive comments of the Editor, Jianqing Fan, an Associate Editor and two referees.

Department of Mathematical Sciences
The University of Texas at Dallas
Richardson, Texas 75083-0688
USA
E-mail: efrom@utdallas.edu